\documentclass[preprint,12pt]{elsarticle}
\usepackage{}

\usepackage{amsmath}
\usepackage{cases}
\usepackage{color}
\usepackage{hyperref}
\usepackage{amssymb}
\usepackage{graphicx}
\usepackage{graphics}
\usepackage{subfigure}
\usepackage{slashbox}
\usepackage{appendix}
\usepackage{textcomp}
\usepackage{bm}
\usepackage{url}

\makeatletter
  \newcommand\figcaption{\def\@captype{figure}\caption}
  \newcommand\tabcaption{\def\@captype{table}\caption}
\makeatother

\biboptions{numbers,sort&compress}

\begin{document}

\begin{center}

{\Large Analytic Solutions of Von K{\'a}rm{\'a}n Plate \\ under Arbitrary Uniform Pressure \\ 
Part (II): Equations in Integral Form}

\vspace{0.3cm}

Xiaoxu Zhong $^3$, Shijun Liao $^{1,2,3}$  \footnote{Corresponding author.  Email address: sjliao@sjtu.edu.cn}

\vspace{0.3cm}

$^1$ State Key Laboratory of Ocean Engineering, Shanghai 200240, China\\
$^2$ Collaborative Innovative Center for Advanced Ship and Deep-Sea Exploration, Shanghai 200240, China\\
$^3$ School of Naval Architecture, Ocean and Civil Engineering\\
Shanghai Jiao Tong University, Shanghai 200240, China

 \end{center}

\hspace{-0.6cm}{\bf Abstract}  {
\em In this paper, the homotopy analysis method (HAM) is successfully applied to solve the Von K{\'a}rm{\'a}n's plate equations in the integral form for a circular plate with the clamped boundary under an arbitrary uniform external pressure.   Two HAM-based approaches are proposed.  One is for a given external load $Q$, the other for a given central deflection.  Both of them are valid for an  arbitrary uniform external pressure by means of choosing a proper value of the so-called convergence-control parameters $c_1$ and $c_2$ in the frame of the HAM.  Besides,  it is found that iteration can greatly accelerate the convergence of solution series.   In addition,  we prove that the interpolation iterative method \cite{Keller, Zheng} is a special case of the HAM-based $1$st-order iteration approach for a given external load $Q$ when $c_{1}=-\theta$ and $c_{2}=-1$, where $\theta$ denotes the interpolation parameter of the   interpolation iterative method.   Therefore,  like Zheng and Zhou \cite{Zheng2},  one can similarly prove  that  the HAM-based approaches are valid for an arbitrary uniform external pressure, at least in some special cases such as $c_{1}=-\theta$ and $c_{2}=-1$.    Furthermore,  it is found that  the HAM-based iteration approaches converge much faster than the interpolation iterative method \cite{Keller, Zheng}.  All of these illustrate the validity and potential of the HAM for the famous Von K{\'a}rm{\'a}n's plate equations, and show the superiority of the HAM over perturbation methods.    
}

\vspace{0.3cm}

\hspace{-0.6cm}{\bf Key Words} circular plate, \sep uniform external pressure, \sep homotopy analysis method

%% main text
\section{Introduction}

In Part (\uppercase\expandafter{\romannumeral1}), we solved the Von K{\'a}rm{\'a}n's plate equations \cite{Karman} in the differential form for a circular plate under an arbitrary uniform external pressure by means of the homotopy analysis method (HAM)  \cite{liaoPhd, Liaobook, liaobook2, Liao2009, Liao2010, Liao2015MDDiM,KV2012}, an analytical approximation method for highly nonlinear problems.  By means of the so-called convergence-control parameter $c_{0}$, convergent analytic approximations for four kinds of boundaries were obtained, with large enough ratio of central deflection to thickness $w(0)/h>20$.  It is found that the convergence-control parameter $c_0$ plays an important role:  it is the convergence-control parameter $c_{0}$ that can guarantee the convergence of solution series and iteration,  and thus distinguishes the HAM from other analytic methods.  Besides,   it is found that the iteration technique can greatly increase the computational efficiency.  In addition, it  was proved in Part (I) that the perturbation methods for {\em any} a perturbation quantity (including Vincent's \cite{Vincent} and Chien's \cite{Qian} perturbation methods) and the modified iteration method \cite{YehLiu} are only the special cases of the HAM when $c_{0}=-1$.

The Von K{\'a}rm{\'a}n's plate equations \cite{Karman} in the integral form describing the large deflection of a circular thin plate under a uniform external pressure read
\begin{eqnarray}
    \varphi(y) &=& -\int_{0}^{1}K(y,\varepsilon)\cdot \frac{1}{\varepsilon^{2}} \cdot S(\varepsilon)\varphi(\varepsilon)d\varepsilon-\int_{0}^{1}K(y,\varepsilon)\cdot Qd\varepsilon, \label{integral1:origin:U}\\
    S(y) &=& \frac{1}{2}\int_{0}^{1}G(y,\varepsilon)\cdot \frac{1}{\varepsilon^{2}} \cdot \varphi^2(\varepsilon)d\varepsilon, \label{integral2:origin:U}
\end{eqnarray}
in which
\begin{equation}
    K(y,\varepsilon)=\left\{
    \begin{array}{ll}
    (\lambda-1)y\varepsilon+y,  \;\;y\leq\varepsilon,\\
    (\lambda-1)y\varepsilon+\varepsilon,  \;\;y>\varepsilon, \end{array}\right.\label{def:K}
\end{equation}
\begin{equation}
    G(y,\varepsilon)=\left\{
    \begin{array}{ll}
    (\mu-1)y\varepsilon+y,  \;\;y\leq\varepsilon,\\
    (\mu-1)y\varepsilon+\varepsilon,  \;\;y>\varepsilon, \end{array}\right.\label{def:G}
\end{equation}
with the definitions
\begin{equation}
    y=\frac{r^{2}}{R_{a}^{2}},\;\;\;\;W(y)=\sqrt{3(1-\nu^{2})}\frac{w(y)}{h},\;\;\;\;\varphi(y)=y\frac{dW(y)}{dy},  \label{dimensionless1}
\end{equation}
\begin{equation}
    S(y)=3(1-\nu^{2})\frac{R_{a}^{2}N_{r}}{Eh^{3}}y,\;\;\;\;Q=\frac{3(1-\nu^{2})\sqrt{3(1-\nu^{2})}R_{a}^{4}}{4Eh^{4}}p, \label{dimensionless2}
\end{equation}
where $r$ is the radial coordinate whose origin locates at the center of the plate, $E$, $\nu$, $R_{a}$, $h$, $w(y)$, $N_{r}$ and $p$ are Young's modulus of elasticity, the Poisson's ratio, radius, thickness, deflection, the radial membrane force of the plate and the external uniform load, respectively, $\lambda$ and $\mu$ are parameters related to the boundary conditions at $y=1$, $Q$ is a constant related to the uniform external load, respectively.  The dimensionless central deflection
\begin{equation}
    W(y)=-\int_{y}^{1}\frac{1}{z}\varphi(z)d z  \label{deflection}
\end{equation}
can be derived from Eq.(\ref{dimensionless1}). Four kinds of boundaries are considered:

(a)~~Clamped: $\lambda=0$ and $\mu=2/(1-\nu)$;

(b)~~Moveable clamped: $\lambda=0$ and $\mu=0$;

(c)~~Simple support: $\lambda=2/(1+\nu)$ and $\mu=0$;

(d)~~Simple hinged support: $\lambda=2/(1+\nu)$ and $\mu=2/(1-\nu)$.

Keller and Reiss \cite{Keller}  proposed the interpolation iterative method to solve the Von K{\'a}rm{\'a}n's plate equations in the integral form by introducing an interpolation parameter to the iteration procedure, and they  successfully  obtained convergent solutions for loads as high as $Q=7000$.   Further, Zheng and Zhou \cite{Zheng, Zheng2} proved that convergent solutions can be obtained by the interpolation iterative method for an arbitrary uniform external pressure.  Their excellent work is a milestone in this field.

 Iterative procedures of the interpolation iterative method \cite{Keller, Zheng} for the Von K{\'a}rm{\'a}n's plate equations in the integral form   are:
\begin{eqnarray}
\psi_{n}(y) &=& \frac{1}{2}\int_{0}^{1}G(y,\varepsilon)\cdot \frac{1}{\varepsilon^{2}} \cdot\vartheta_{n}^{2}(\varepsilon)d\varepsilon,
\label{AppB:Interpolation:origin1}\\
\vartheta_{n+1}(y) &=& (1-\theta)\vartheta_{n}(y)-\theta\int_{0}^{1}K(y,\varepsilon)Qd\varepsilon\\ \nonumber
    &-&\theta\int_{0}^{1}K(y,\varepsilon)\cdot \frac{1}{\varepsilon^{2}} \cdot\vartheta_{n}(\varepsilon)\psi_{n}(\varepsilon)d\varepsilon.
 \label{AppB:Interpolation:origin2}
\end{eqnarray}
with the definition of the initial guess
\begin{equation}
\vartheta_{1}(y)=-\frac{Q\theta}{2}[(\lambda+1)y-y^{2}]. \label{AppB:Interpolation:initial}
\end{equation}

The homotopy analysis method (HAM) \cite{liaoPhd, Liaobook, liaobook2, Liao2009, Liao2010, Liao2015MDDiM,KV2012} was proposed by Liao \cite{liaoPhd} in $1992$.  Unlike perturbation technique, the HAM is independent of any small/large physical parameters. Besides, unlike other analytic techniques, the HAM provides a simple way to guarantee the convergence of  solution series by means of introducing the so-called ``convergence-control parameter".   In addition, the HAM provides us great freedom to choose equation-type and solution expression of the high-order linear equations.   As a powerful technique to solve highly nonlinear equations, the HAM was successfully employed to solve various types of nonlinear problems over the past two decades \cite{Abbasbandy2006,Hayat2006, KV2008, Liang2010, Ghotbi2011, Nassar2011, Mastroberardino2011, Aureli2014, Duarte2015, Nagarajaiah2015, Gorder2015}.  Especially, as shown in \cite{xu2012JFM, Liu2014JFM, Liu2015JFM, Liao2016JFM, BookChap3-Liao2015}, the HAM  can bring us something completely new/different: the steady-state resonant waves were first predicted by the HAM in theory and then confirmed experimentally in a lab \cite{Liu2015JFM}. 

In this paper, we propose two approaches in the frame of the HAM to solve the Von K{\'a}rm{\'a}n's plate equations in the integral form.  Besides, we proof that   the interpolation iterative method \cite{Keller, Zheng} is only a special case of the HAM.   

\section{HAM approach for given external load $Q$}

\subsection{Mathematical formulas}

Following Zheng \cite{Zheng}, we express $\varphi(y)$ and $S(y)$ as
\begin{equation}
   \varphi(y)=\sum_{k=1}^{+\infty} a_{k}\cdot y^{k},\;\;\;\;S(y)=\sum_{k=1}^{+\infty}b_{k}\cdot y^{k}, \label{homotopy:series}
\end{equation}
where $a_{k}$ and $b_{k}$ are constant coefficients to be determined. They provide us the so-called  ``solution expression" of $\varphi(y)$ and $S(y)$ in the frame of the HAM.

Let $\varphi_{0}(y)$ and $S_{0}(y)$ be initial guesses of $\varphi(y)$ and $S(y)$. Moreover, let $c_{1}$ and $c_{2}$ denote the non-zero auxiliary parameters, called the convergence-control parameters, and $q\in[0,1]$ the embedding parameter, respectively.  Then we construct a family of differential equations in $q\in[0,1]$, namely the zeroth-order deformation equation:
\begin{eqnarray}
   (1-q)[\Phi(y;q)-\varphi_{0}(y)] &=& c_{1} \; q \; {\cal N}_{1}(y;q), \label{AppB:HAM:zeroth1}\\
   (1-q)[\Xi(y;q)-S_{0}(y)] &=& c_{2} \; q \; {\cal N}_{2}(y;q), \label{AppB:HAM:zeroth2}
\end{eqnarray}
where
\begin{eqnarray}
    {\cal N}_{1}(y;q) &=& \Phi(y;q)+\int_{0}^{1}K(y,\varepsilon)\cdot \frac{1}{\varepsilon^{2}} \cdot \Phi(\varepsilon ;q)\Xi(\varepsilon ;q)d\varepsilon \nonumber \\
    &+& \int_{0}^{1}K(y,\varepsilon)\cdot Qd\varepsilon, \label{AppB:HAM:operator1}\\
    {\cal N}_{2}(y;q) &=& \Xi(y;q)-\frac{1}{2}\int_{0}^{1}G(y,\varepsilon)\cdot \frac{1}{\varepsilon^{2}} \cdot \Phi^2(\varepsilon ;q)d\varepsilon, \label{AppB:HAM:operator2}
\end{eqnarray}
are two nonlinear operators.

When $q=0$, Eqs.~(\ref{AppB:HAM:zeroth1}) and (\ref{AppB:HAM:zeroth2}) have the solution
\begin{equation}
    \Phi(y;0)=\varphi_{0}(y),\;\;\;\;\Xi(y;0)=S_{0}(y).
\end{equation}
 When $q=1$, they are equivalent to the original equations (\ref{integral1:origin:U})-(\ref{integral2:origin:U}), provided
\begin{equation}
    \Phi(y;1)=\varphi(y),\;\;\;\;\Xi(y;1)=S(y). \label{embedding:parameter}
\end{equation}
Therefore, as $q$ increases from $0$ to $1$, $\Phi(y;q)$ varies continuously from the initial guess $\varphi_{0}(y)$ to the solution $\varphi(y)$, so does $\Xi(y;q)$ from the initial guess $S_{0}(y)$ to the solution $S(y)$. In topology, these kinds of continuous variation are called deformation.  That is why Eqs.~(\ref{AppB:HAM:zeroth1})-(\ref{AppB:HAM:zeroth2}) constructing the homotopies $\Phi(y;q)$ and $\Xi(y;q)$ are called the zeroth$-$order deformation equations.

Expanding $\Phi(y;q)$ and $\Xi(y;q;a)$ into Taylor series with respect to the embedding parameter $q$, we have the so-called homotopy-series:
\begin{equation}
    \Phi(y;q)=\varphi_{0}(y)+\sum_{k=1}^{+\infty}\varphi_{k}(y)q^{k}, \label{socalled:homotopy:phi}
\end{equation}
\begin{equation}
    \Xi(y;q)=S_{0}(y)+\sum_{k=1}^{+\infty}S_{k}(y)q^{k}, \label{socalled:homotopy:s}
\end{equation}
where
\begin{equation}
   \varphi_{k}(y)={\cal D}_{k}[\Phi(y;q)], \;\;\; S_{k}(y)={\cal D}_{k}[\Xi(y;q)], \label{homotopy:S}
\end{equation}
in which 
\begin{equation}
   {\cal D}_{k}[f]=\frac{1}{k!}\frac{\partial^{k}f}{\partial q^{k}}\bigg{|}_{q=0} \label{Dm}
\end{equation}
is called the $k$th-order homotopy-derivative of $f$.

Note that there are two convergence-control parameters $c_{1}$ and $c_{2}$  in the homotopy-series (\ref{socalled:homotopy:phi}) and (\ref{socalled:homotopy:s}), respectively. Assume that  $c_{1}$ and $c_{2}$ are properly chosen so that the homotopy-series (\ref{socalled:homotopy:phi}) and (\ref{socalled:homotopy:s}) converge at $q=1$. Then according to Eq.~(\ref{embedding:parameter}), the so-called homotopy-series solutions read:
\begin{equation}
    \varphi(y)=\varphi_{0}(y)+\sum_{k=1}^{+\infty}\varphi_{k}(y),
\end{equation}
\begin{equation}
    S(y)=S_{0}(y)+\sum_{k=1}^{+\infty}S_{k}(y).
\end{equation}
 The $n$th-order approximation of $\varphi(y)$ and $S(y)$ read
\begin{equation}
   \tilde{\Phi}_{n}(y)=\sum_{k=0}^{n}\varphi_{k}(y),\label{sum:phi:n}
\end{equation}
\begin{equation}
   \tilde{\Xi}_{n}(y)=\sum_{k=0}^{n}S_{k}(y).\label{sum:s:n}
\end{equation}

Substituting the power series (\ref{socalled:homotopy:phi}) and (\ref{socalled:homotopy:s}) into the zeroth-order deformation equations (\ref{AppB:HAM:zeroth1}) and (\ref{AppB:HAM:zeroth2}), and then equating the like-power of the embedding parameter $q$, we have the so-called $k$th-order deformation equations
\begin{equation}
    \varphi_{k}(y)-\chi_{k}\varphi_{k-1}(y)=c_{1}\; \delta_{1,k-1}(y), \label{AppB:HAM:highorder1}
\end{equation}
\begin{equation}
    S_{k}(y)-\chi_{k}S_{k-1}(y)=c_{2} \; \delta_{2,k-1}(y), \label{AppB:HAM:highorder2}
\end{equation}
where
\begin{eqnarray}
    \delta_{1,k-1}(y) &=& {\cal D}_{k-1}[{\cal N}_{1}(y;q)]\nonumber \\
    &=& \varphi_{k-1}(y)+\int_{0}^{1}K(y,\varepsilon)\cdot \frac{1}{\varepsilon^{2}} \cdot \sum_{i=0}^{k-1}\varphi_{i}(\varepsilon)S_{k-1-i}(\varepsilon)d\varepsilon \nonumber\\
    &+& (1-\chi_{k})\int_{0}^{1}K(y,\varepsilon)\cdot Qd\varepsilon, \label{AppB:HAM:delta1}\\
    \delta_{2,k-1}(y)&=& {\cal D}_{k-1}[{\cal N}_{2}(y;q)] \nonumber \\
    &=& S_{k-1}(y)-\frac{1}{2}\int_{0}^{1}G(y,\varepsilon)\cdot \frac{1}{\varepsilon^{2}} \cdot\sum_{i=0}^{k-1} \varphi_{i}(\varepsilon)\varphi_{k-1-i}(\varepsilon)d\varepsilon, \label{AppB:HAM:delta2}
\end{eqnarray}
with the definition 
\begin{equation}
\chi_m =\left\{
\begin{array}{cc}
1 & \mbox{when $m=1$}, \\
0 & \mbox{when $m \geq 2$.}
\end{array}
\right. \label{def:chi}
\end{equation}
Note that the $k$th-order deformation equations (\ref{AppB:HAM:highorder1})-(\ref{AppB:HAM:highorder2}) are linear.

Note that  we have great freedom to choose the initial guesses $\varphi_{0}(y)$ and $S_{0}(y)$. According to \cite{Liaobook}, $\varphi_{0}(y)$ and $S_{0}(y)$ should obey the solution expression (\ref{homotopy:series}), thus, we choose the initial guesses:
\begin{eqnarray}
\varphi_{0}(y) &=& \frac{Q c_{0}}{2}[(\lambda+1)y-y^{2}],\\
S_{0}(y) &=& 0,
\end{eqnarray}
which satisfy all boundary conditions.   Besides, for the sake of simplicity, set 
\begin{equation}
c_{1}=c_{2}=c_{0}.
\end{equation}
In order to measure the accuracy of approximations,  we define the sum of the two discrete squared residuals
\begin{equation}
   Err=\frac{1}{K+1}\sum_{i=0}^{K}\left\{\left[{\cal N}_{1}\left(\frac{i}{K}\right)\right]^{2}+\left[{\cal N}_{2}\left(\frac{i}{K}\right)\right]^{2}\right\}, \label{discrete:residual}
\end{equation}
where ${\cal N}_1$ and ${\cal N}_2$ are two nonlinear operators defined by (\ref{AppB:HAM:operator1}) and (\ref{AppB:HAM:operator2}), respectively, and  $K=100$ is used in the whole of this paper.  Obviously, the smaller the $Err$, the more accurate the approximation.

\subsection{Convergent results given by the HAM approach without iteration}

Without loss of generality, the Von K{\'a}rm{\'a}n's plate equations with the clamped boundary is studied at first, and the Poisson's ratio $\nu$ is taken to be 0.3 in all cases considered in this paper.

\begin{figure}
    \begin{center}
        \begin{tabular}{cc}
            \includegraphics[width=3in]{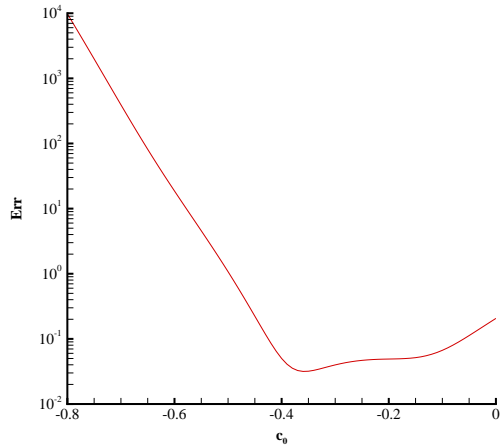} \\
        \end{tabular}
    \caption{The sum of the squared residual errors $Err$ versus $c_{0}$ in the case of the given external load $Q = 5$ for a circular plate with the clamped boundary, given by the HAM-based approach (without iteration).} \label{HAM:Q5:opt:c0:5th}
    \end{center}
\end{figure}

First of all, the HAM-based approach (without iteration) is used to solve the Von K{\'a}rm{\'a}n's plate equations for a given external load $Q$ with the clamped boundary. Take the case of $Q=5$ as an example. The optimal value of $c_{0}$ is determined by the minimum of $Err$ defined by (\ref{discrete:residual}), i.e. the sum of the squared residual errors of the two governing equations. As shown in Fig.~\ref{HAM:Q5:opt:c0:5th}, the sum of the squared residual errors arrives its minimum at $c_{0}\approx-0.35$, which suggests that the optimal value of $c_{0}$ is about $-0.35$.  As shown in Table~\ref{Table:HAM:Q5}, the sum of the squared residual errors quickly decreases to $1.7\times10^{-7}$  by means of $c_{0}\approx-0.35$ in the case of $Q=5$, with the convergent ratio of the maximum central deflection to plate thickness  $w(0)/h=0.62$.  Note that Vincent's perturbation result \cite{Vincent} (using $Q$ as the perturbation quantity) for a circular plate with the clamped boundary is only valid for $w(0)/h<0.52$, corresponding to $Q<3.9$, and thus fails in the case of $Q=5$.   So, the convergence control parameter $c_0$  indeed provides us a simple way to guarantee the convergence of solution series.  This illustrates that our HAM-based approach (mentioned in \S~2.1)  has advantages over the perturbation method \cite{Vincent}.

For given values of $Q$, we can always find its optimal convergence control parameter $c_{0}$ in a similar way, which can be expressed by the empirical formula:
\begin{equation}
c_{0}=-\frac{13}{13+Q^2}. \label{c0:noniter:Q}
\end{equation}
Besides, the convergent homotopy-approximations of central deflection $w(0)/h$ in case of different values of $Q$ for a circular plate with the clamped boundary are given in Table {\ref{HAM:noniter:diff:Q}}.

\begin{table}[t]
\tabcolsep 0pt
\caption{The sum of the squared residual errors $Err$ and the central deflection $w(0)/h$ versus the order of approximation in the case of  $Q = 5$ for a circular plate with the clamped boundary, given by the HAM approach (see \S~2.1) without iteration using $c_{0} = -0.35$.}
\vspace*{-12pt}\label{Table:HAM:Q5}
\begin{center}
\def\temptablewidth{0.8\textwidth}
{\rule{\temptablewidth}{1pt}}
\begin{tabular*}{\temptablewidth}{@{\extracolsep{\fill}}cccccc}
$m$, order of approx. &  $Err$    &  $w(0)/h$                    \\
\hline
10 &  $3.3\times10^{-4}$  &  0.64                          \\
20 &  $6.5\times10^{-5}$  &  0.62                  \\
30 &  $1.1\times10^{-5}$ &  0.62                  \\
40 &  $1.6\times10^{-6}$ &  0.62                         \\
50 & $1.7\times10^{-7}$   &  0.62                      \\
\end{tabular*}
{\rule{\temptablewidth}{1pt}}
 \end{center}
 \end{table}

\begin{table}
\tabcolsep 0pt
\caption{The convergent results of the central deflection $w(0)/h$ in the case of the different values of $Q$ for a circular plate with the clamped boundary, given by the HAM-based approach (see \S~2.1) without iteration  using the optimal convergence-control parameter $c_{0}$ given by (\ref{c0:noniter:Q}).}
\vspace*{-12pt}\label{HAM:noniter:diff:Q}
\begin{center}
\def\temptablewidth{0.8\textwidth}
{\rule{\temptablewidth}{1pt}}
\begin{tabular*}{\temptablewidth}{@{\extracolsep{\fill}}cccccc}
$Q$  & $c_{0}$ &  $w(0)/h$       \\
\hline
1    & -0.93   &  0.15        \\
2   & -0.76   &  0.29      \\
3   & -0.59   &  0.41     \\
4   & -0.45  &  0.53     \\
5   &-0.34   &  0.62   \\
\end{tabular*}
{\rule{\temptablewidth}{1pt}}
 \end{center}
 \end{table}

\subsection{Convergence acceleration by iteration}

According to Liao \cite{Liaobook}, the convergence of the homotopy-series solutions can be greatly accelerated by means of iteration technique, which uses the $M$th-order homotopy-approximations
\begin{eqnarray}
  \varphi^*(y) &\approx& \varphi_{0}(y)+\sum_{k=1}^{M}\varphi_{k}(y), \label{Q:M:order:phi}\\
  S^*(y) &\approx& S_{0}(y)+\sum_{k=1}^{M}S_{k}(y), \label{Q:M:order:s}
\end{eqnarray}
as the new initial guesses of $\varphi_{0}(y)$ and $S_{0}(y)$ for the next iteration.  This provides us the $M$th-order iteration approach of the HAM.  

Note that the length of the solution expressions increases exponentially in iteration.   To avoid this,  we truncate the right-hand sides  of  Eqs.(\ref{AppB:HAM:highorder1}) and (\ref{AppB:HAM:highorder2}) in the following way 
\begin{equation}
\delta_{1,k}(y)\approx \sum_{m=0}^{N}A_{k,m}\cdot y^{k}, \;\;\;  
\delta_{2,k}(y;a)\approx \sum_{m=0}^{N}B_{k,m}\cdot y^{k}, \label{Q:nterms2}
\end{equation}
where $A_{k,m}$ and $B_{k,m}$ are constant coefficients, and $N$ is called the truncation order, respectively.

\begin{figure}[t]
    \begin{center}
        \begin{tabular}{cc}
            \includegraphics[width=2.5in]{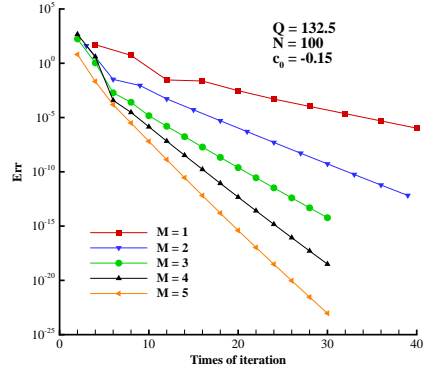} \\
        \end{tabular}
    \caption{The sum of the squared residual errors $Err$ versus the  times of iteration in the case of $Q=132.2$ for a circular plate with the clamped boundary, given by the HAM-based iteration approach (with different orders of iteration) using  the convergence-control parameter $c_{0}=-0.15$ and the truncation order $N=100$. Square: 1st-order; Triangle down: 2nd-order; Circle: 3rd-order; Triangle up: 4th-order; Triangle left: 5th-order.} \label{Iter:Q:order}
    \end{center}
\end{figure}

\begin{figure}
    \begin{center}
        \begin{tabular}{cc}
            \includegraphics[width=2.5in]{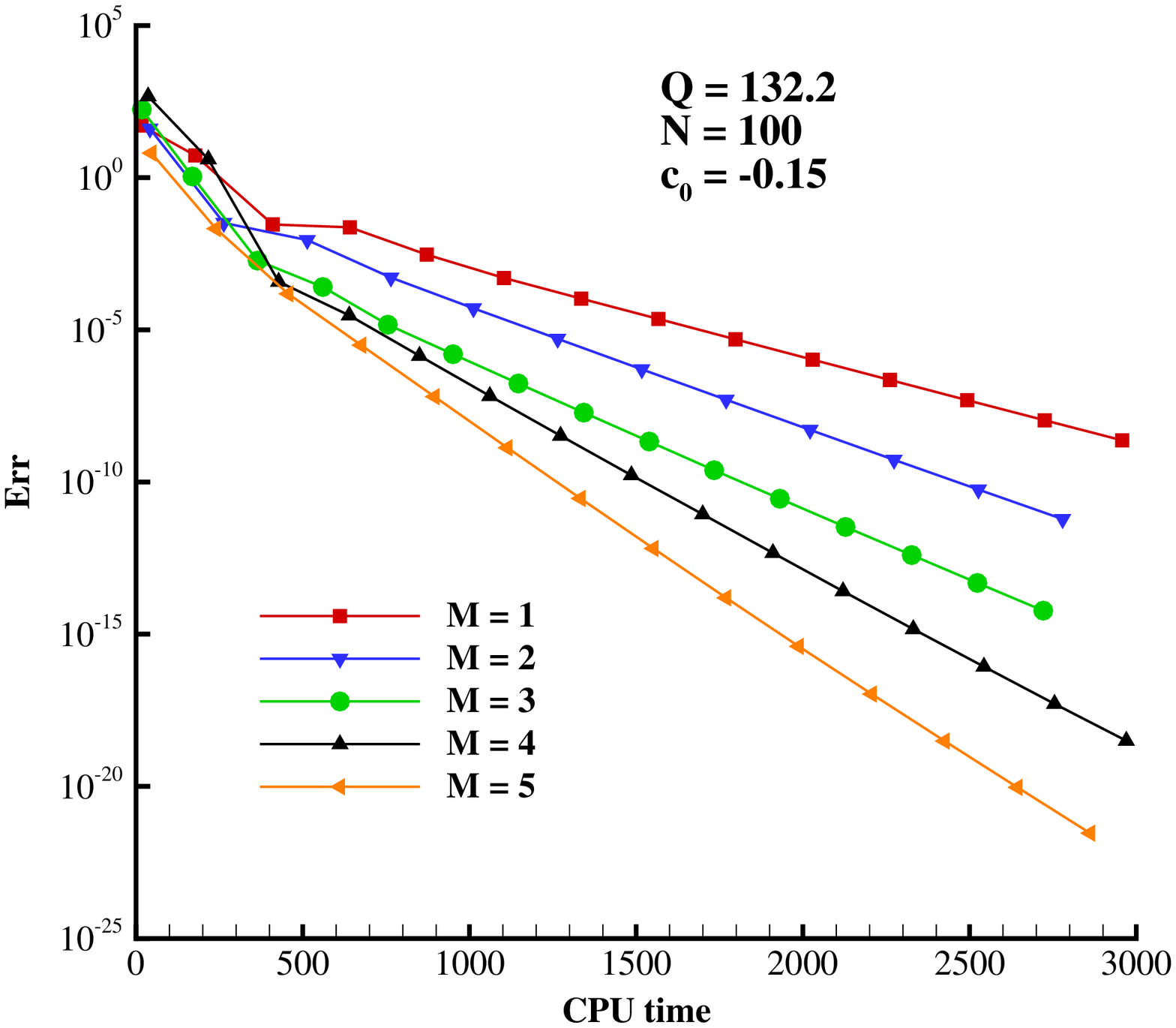} \\
        \end{tabular}
    \caption{The sum of the squared residual errors $Err$ versus the CPU time in the case of $Q=132.2$ for a circular plate with the clamped boundary, given by the HAM-based iteration approach (with different orders of iteration) using  the convergence-control parameter $c_{0}=-0.15$ and the truncation order $N=100$. Square: 1st-order; Triangle down: 2nd-order; Circle: 3rd-order; Triangle up: 4th-order; Triangle left: 5th-order.} \label{Iter:Q:CPU}
    \end{center}
\end{figure}

 Without loss of generality, let us consider the case of $Q=132.2$, corresponding to $w(0)/h=3.0$.  As shown in Fig.~\ref{Iter:Q:order}, the higher the order $M$ of iteration, the less  times of iteration is required for a given accuracy-level of approximation.   Besides, the sum of the squared residual errors $Err$ versus the CPU time for different order of iteration ($M$) are given in Fig.~\ref{Iter:Q:CPU}.  Note that the higher the order of iteration, the faster the approximation converges.  Thus, it is natural for us to choose the 5th-order  iteration approach (i.e. $M=5$) from the view-point of computational efficiency.
 
\begin{figure}
    \begin{center}
        \begin{tabular}{cc}
            \includegraphics[width=2.5in]{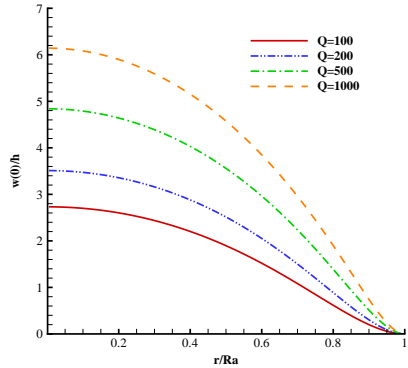} \\
        \end{tabular}
    \caption{The deflection curves of a circular plate with the clamped boundary given by the HAM-based 5th-order iteration approach in the case of $Q=100, 200, 500, 1000$. Solid line: $Q=100$; Dash-double-dotted line: $Q=200$; Dash-dotted line: $Q=500$; Dashed line: $Q=1000$.} \label{HAM:Q:deflection}
    \end{center}
\end{figure}

For a given $Q$,  we use the truncation order $N=100$.  It is found that the corresponding optimal convergence-control parameter $c_{0}$ can be expressed by the empirical formula:
\begin{equation}
  c_{0}=-\frac{23}{Q+23}. \label{c0:iter:Q}
\end{equation}
As shown in Table~\ref{HAM:iter:Q1000}, convergent analytic solutions can be obtained even in the case of $Q=1000$, corresponding to $w(0)/h=6.1$.   Here, it should be emphasized that Vincent's perturbation result \cite{Vincent} (using $Q$ as the perturbation quantity) for a circular plate with the clamped boundary is only valid for $w(0)/h<0.52$, corresponding to $Q<3.9$ only.  This illustrates again the superiority of the HAM over perturbation methods.    The convergent approximations of central deflection $w(0)/h$ in case of different values of $Q$ for a circular plate with the clamped boundary are listed in Table \ref{HAM:iter:diff:Q}.  In addition, several deflection curves under different loads $Q$ are shown in Fig.~\ref{HAM:Q:deflection}.  

\begin{table}[t]
\tabcolsep 0pt
\caption{The sum of the squared residual errors $Err$ and the homotopy-approximations of $w(0)/h$ at the $m$th iteration  in the case of $Q=1000$ for a circular plate with the clamped boundary, given by the HAM-based 5th-order iteration approach with the optimal convergence-control parameter $c_{0}=-0.02$ and the truncation order $N=100$.}
\vspace*{-12pt}\label{HAM:iter:Q1000}
\begin{center}
\def\temptablewidth{0.8\textwidth}
{\rule{\temptablewidth}{1pt}}
\begin{tabular*}{\temptablewidth}{@{\extracolsep{\fill}}cccccc}
$m$, times of iteration  & $Err$ &  $w(0)/h$   \\
\hline
20    & $2.0\times10^{-1}$  &  6.1        \\
40   & $1.4\times10^{-3}$   &  6.1      \\
60   & $1.4\times10^{-5}$   &  6.1     \\
80   & $1.4\times10^{-7}$  &   6.1     \\
100   &$1.5\times10^{-9}$   &  6.1   \\
\end{tabular*}
{\rule{\temptablewidth}{1pt}}
 \end{center}
 \end{table}

\begin{table}
\tabcolsep 0pt
\caption{The convergent homotopy-approximation of the central deflection $w(0)/h$ in case of different values of $Q$ for a circular plate with the clamped boundary, given by the HAM-based 5th-order iteration approach using the optimal convergence-control parameter $c_{0}$ given in (\ref{c0:iter:Q}) and the truncation order $N=100$.}
\vspace*{-12pt}\label{HAM:iter:diff:Q}
\begin{center}
\def\temptablewidth{0.8\textwidth}
{\rule{\temptablewidth}{1pt}}
\begin{tabular*}{\temptablewidth}{@{\extracolsep{\fill}}cccccc}
$Q$  & $c_{0}$ &  $w(0)/h$       \\
\hline
200    & -0.10   &  3.5        \\
400   & -0.05   &  4.5      \\
600   & -0.04   &  5.2     \\
800   & -0.03  &   5.7     \\
1000   &-0.02   &  6.1   \\
\end{tabular*}
{\rule{\temptablewidth}{1pt}}
 \end{center}
 \end{table}

\subsection{Relations between the HAM and the interpolation iterative method}

Here we prove that the interpolation iterative method \cite{Keller, Zheng}  is  a special case of the HAM-based 1st-order iteration approach.

Let $\theta$ denote the interpolation parameter in the interpolation iterative method \cite{Keller, Zheng}.    Set $c_{1}=-\theta$ and $c_{2}=-1$, since we have the great freedom to choose their values in the frame of the HAM.   According to (\ref{AppB:HAM:highorder1})-(\ref{def:chi}) and (\ref{Q:M:order:phi})-(\ref{Q:M:order:s}), we have the $1$st-order homotopy-approximations of $\varphi(y)$ and $S(y)$:
\begin{eqnarray}
    \varphi^*(y)&=&\varphi_{0}(y)+\varphi_{1}(y) \nonumber \\
    &=&\varphi_{0}(y)-\theta\delta_{1,1}(y) \nonumber \\
    &=&(1-\theta)\varphi_{0}(y)-\theta\int_{0}^{1}K(y,\varepsilon)Qd\varepsilon \nonumber \\
    &-&\theta\int_{0}^{1}K(y,\varepsilon)\cdot \frac{1}{\varepsilon^{2}} \cdot\varphi_{0}(\varepsilon)S_{0}(\varepsilon)d\varepsilon, 
 \label{AppB:HAM:SumPhinew1}
\end{eqnarray}
\begin{eqnarray}
S^*(y)&=&S_{0}(y)+S_{1}(y) \nonumber \\
&=&S_{0}(y)-\delta_{2,1}(y) \nonumber \\
&=&\frac{1}{2}\int_{0}^{1}G(y,\varepsilon)\cdot \frac{1}{\varepsilon^{2}} \cdot\varphi_{0}^{2}(\varepsilon)d\varepsilon.
\label{AppB:HAM:SumSnew1}
\end{eqnarray}

Assume that $\varphi_0(y)$ and $S_0(y)$ are known.  We use the following iterative procedure:
\begin{enumerate}
\item[(A)]   Calculate $S^*(y)$ according to Eq.~(\ref{AppB:HAM:SumSnew1});
\item[(B)]    $S_{0}(y)$ is replaced by $S^*(y)$ that is used as the new initial guess;
\item[(C)]    Calculate $\varphi^*(y)$ according to Eq.~(\ref{AppB:HAM:SumPhinew1});
\item[(D)]  $\varphi_0(y)$ is replaced by $\varphi^*(y)$ that is used as the new initial guess.  
\end{enumerate}
At the $n$th times of  iteration, write
\begin{equation}
  \hat{\Phi}_{n}(y)=\varphi^*(y),\;\;\;\;\;\;\;    \hat{\Xi}_{n-1}(y)=S^*(y).\nonumber
\end{equation}
Then,  the procedure of the HAM-based 1st-order iteration approach is expressed by
\begin{eqnarray}
\hat{\Xi}_{n-1}(y) &=& \frac{1}{2}\int_{0}^{1}G(y,\varepsilon)\cdot \frac{1}{\varepsilon^{2}} \cdot\hat{\Phi}_{n-1}^{2}(\varepsilon)d\varepsilon,
\label{AppB:HAM:origin1}\\
\hat{\Phi}_{n}(y) &=& (1-\theta)\hat{\Phi}_{n-1}(y)-\theta\int_{0}^{1}K(y,\varepsilon)Qd\varepsilon\nonumber \\
    &-& \theta\int_{0}^{1}K(y,\varepsilon)\cdot \frac{1}{\varepsilon^{2}} \cdot\hat{\Phi}_{n-1}(\varepsilon)\hat{\Xi}_{n-1}(\varepsilon)d\varepsilon.
\label{AppB:HAM:origin2}
\end{eqnarray}
We choose the initial guess
\begin{equation}
\hat{\Phi}_{0}(y)=-\frac{Q\theta}{2}[(\lambda+1)y-y^{2}].\label{AppB:HAM:initial}
\end{equation}

Note that, the iterative procedures (\ref{AppB:Interpolation:origin1}), (\ref{AppB:Interpolation:origin2}) and the initial solution (\ref{AppB:Interpolation:initial}) of the interpolation iterative method \cite{Keller, Zheng} are exactly the same as those of the HAM-based iteration approach (\ref{AppB:HAM:origin1})-(\ref{AppB:HAM:initial}).  Thus, the interpolation iterative method \cite{Keller, Zheng} is a special case of the HAM-based $1$st-order iteration approach when $c_{1}=-\theta$ and $c_{2}=-1$. It should be emphasized that, as shown in Figs.~\ref{Iter:Q:order} and \ref{Iter:Q:CPU}, the higher the order of the iteration, the faster the approximations converge. So, the interpolation iterative method \cite{Keller, Zheng}  should correspond to the  slowest  one  among  the  HAM-based $M$th-order iteration approaches (up to $M=5$).   Finally, it should be emphasized that the interpolation iterative method \cite{Keller, Zheng} is valid for an {\em arbitrary} uniform external pressure, as proved by Zheng and Zhou \cite{Zheng2}. So, following Zheng and Zhou \cite{Zheng2},  one could prove that the HAM-based approach mentioned in \S~2.1 is valid for an arbitrary uniform pressure at least in some special cases, such as $c_{1}=-\theta$ and $c_{2}=-1$.  This reveals the important role of the convergence-control parameters $c_1$ and $c_2$ in the frame of the HAM.   

\section{HAM approach for given central deflection}

According to Chien \cite{Qian},  it makes sense to introduce the central deflection into the Von K{\'a}rm{\'a}n's plate equations for a circular plate  so as to enlarge the convergent region. Based on this knowledge, the HAM-based approach for given central deflection is employed to solve the Von K{\'a}rm{\'a}n's plate equations in the integral form with the clamped boundary.

\subsection{Mathematical formulas}

Given
\begin{equation}
    W(0)=a, \nonumber
\end{equation}
we have due to Eq.~(\ref{deflection}) an additional equation for the corresponding unknown value of $Q$:
\begin{equation}
    \int_{0}^{1}\frac{1}{\varepsilon}\varphi(\varepsilon)d\varepsilon=-a. \label{restricted:condition}
\end{equation}
Let $\varphi_{0}(y)$ and $S_{0}(y)$ denote initial guesses of $\varphi(y)$ and $S(y)$,  which satisfy the restriction condition (\ref{restricted:condition}), $c_{1}$ and $c_{2}$ the non-zero auxiliary parameters, called the convergence-control parameters, $q\in[0,1]$ the embedding parameter, respectively. Then,  we construct the so-called zeroth-order deformation equations
\begin{equation}
   (1-q)[\tilde{\Phi}(y;q)-\varphi_{0}(y)]=c_{1}\; q \; \tilde{{\cal N}}_{1}(y;q), \label{zeroth:deformation:one}
\end{equation}
\begin{equation}
   (1-q)[\tilde{\Xi}(y;q)-S_{0}(y)]=c_{2}\; q \; \tilde{{\cal N}}_{2}(y;q), \label{zeroth:deformation:two}
\end{equation}
subject to the restricted condition
\begin{equation}
    \int_{0}^{1}\frac{1}{\varepsilon}\tilde{\Phi}(\varepsilon;q) d\varepsilon=-a. \label{W:restricted:condition}
\end{equation}
where 
\begin{eqnarray}
    \tilde{{\cal N}}_{1}(y;q)&=&\tilde{\Phi}(y;q)+\int_{0}^{1}K(y,\varepsilon)\cdot \frac{1}{\varepsilon^{2}} \cdot \tilde{\Phi}(\varepsilon ;q)\tilde{\Xi}(\varepsilon ;q)d\varepsilon \nonumber \\
    &+&\int_{0}^{1}K(y,\varepsilon)\cdot \tilde{\Theta}(q)d\varepsilon, \label{nonlinear:equation:one}\\
    \tilde{{\cal N}}_{2}(y;q) &=& \tilde{\Xi}(y;q)-\frac{1}{2}\int_{0}^{1}G(y,\varepsilon)\cdot \frac{1}{\varepsilon^{2}} \cdot \tilde{\Phi}^2(\varepsilon ;q)d\varepsilon, \label{nonlinear:equation:two}
\end{eqnarray}
are the two nonlinear operators,  corresponding to  Eqs.~(\ref{integral1:origin:U}) and (\ref{integral2:origin:U}), respectively.  
 Note that  we  introduce here a continuous variation $\tilde{\Theta}(q)$ from the initial guess $Q_0$ to $Q$, since $Q$ is unknown for a given central deflection $a$.    

When $q=0$,  it holds 
\begin{equation}
    \tilde{\Phi}(y;0)=\varphi_{0}(y),\;\;\;\;\tilde{\Xi}(y;0)=S_{0}(y),
\end{equation}
since  $\varphi_{0}(y)$ satisfies the restriction condition.   When $q=1$, they are equivalent to the original equations (\ref{integral1:origin:U}), (\ref{integral2:origin:U}) and (\ref{restricted:condition}), provided
\begin{equation}
    \tilde{\Phi}(y;1)=\varphi(y),\;\;\;\; \tilde{\Xi}(y;1)=S(y),\;\;\;\; \tilde{\Theta}(1)=Q.\label{provide:Q}
\end{equation}
Then, as $q$ increases from 0 to 1, $\tilde{\Phi}(y;q)$ varies continuously from the initial guess $\varphi_{0}(y)$ to the solution $\varphi(y)$, so do $\tilde{\Xi}(y;q)$ from the initial guess $S_{0}(y)$ to the solution $S(y)$, $\tilde{\Theta}(q)$ from the initial guess $Q_{0}$ to the unknown load $Q$ for the given central deflection $a$, respectively.  

Expanding $\tilde{\Phi}(y;q)$, $\tilde{\Xi}(y;q)$ and $\tilde{\Theta}(q)$ into Taylor series with respect to the embedding parameter $q$, we have the so-called homotopy-series:
\begin{eqnarray}
    \tilde{\Phi}(y;q)&=&\varphi_{0}(y)+\sum_{k=1}^{+\infty}\varphi_{k}(y)q^{k}, \label{new:homotopy:phi}\\
    \tilde{\Xi}(y;q)&=&S_{0}(y)+\sum_{k=1}^{+\infty}S_{k}(y)q^{k}, \label{new:homotopy:s}\\
    \tilde{\Theta}(q)&=&Q_{0}+\sum_{k=1}^{+\infty}Q_{k} q^{k}, \label{new:homotopy:Q}
\end{eqnarray}
in which
\begin{equation}
   \varphi_{k}(y)=D_{k}[\Phi(y;q)],\;\;\; S_{k}(y)=D_{k}[\Xi(y;q)],\;\;\; Q_{k}=D_{k}[\Theta(q)], \label{new:Q}
\end{equation}
with the definition ${\cal D}_{k}$ by (\ref{Dm}).  Assume that the homotopy-series (\ref{new:homotopy:phi})-(\ref{new:homotopy:Q}) are convergent at $q=1$. According to (\ref{provide:Q}),  we have the so-called homotopy-series solutions:
\begin{eqnarray}
    \varphi(y)&=&\varphi_{0}(y)+\sum_{k=1}^{+\infty}\varphi_{k}(y),\\
    S(y)&=&S_{0}(y)+\sum_{k=1}^{+\infty}S_{k}(y),\\
    Q&=&Q_{0}+\sum_{k=1}^{+\infty}Q_{k}.
\end{eqnarray}
Substituting the power series (\ref{new:homotopy:phi})-(\ref{new:homotopy:Q}) into the zeroth-order deformation equations (\ref{zeroth:deformation:one})-(\ref{W:restricted:condition}), and then equating the like-power of $q$, we have the so-called $k$th-order deformation equations
\begin{eqnarray}
    \varphi_{k}(y)-\chi_{k}\varphi_{k-1}(y)&=&c_{1} \; \delta_{1,k-1}(y), \label{high:1}\\
    S_{k}(y)-\chi_{k}S_{k-1}(y)&=&c_{2} \; \delta_{2,k-1}(y), \label{high:2}
\end{eqnarray}
subject to the restriction condition
\begin{equation}
\int_{0}^{1}\frac{1}{\varepsilon}\varphi_{k}(\varepsilon)d\varepsilon=0,\label{restrict}
\end{equation}
where $\chi_{k}$ is defined by (\ref{def:chi}), and 
\begin{eqnarray}
    \delta_{1,k-1}(y)&=&{\cal D}_{k-1}[{\cal N}_{1}(y;q)] \nonumber \\
    &=&\varphi_{k-1}(y)+\sum_{i=0}^{k-1} \int_{0}^{1}K(y,\varepsilon)\cdot \frac{1}{\varepsilon^{2}} \cdot \varphi_{i}(\varepsilon)S_{k-1-i}(\varepsilon)d\varepsilon \nonumber \\
    &+&\int_{0}^{1}K(y,\varepsilon)\cdot Q_{k-1}d\varepsilon, \label{AppB:HAM:newdelta1}\\
    \delta_{2,k-1}(y)&=&{\cal D}_{k-1}[{\cal N}_{2}(y;q)] \nonumber \\
    &=&S_{k-1}(y)-\sum_{i=0}^{k-1} \frac{1}{2}\int_{0}^{1}G(y,\varepsilon)\cdot \frac{1}{\varepsilon^{2}} \cdot\varphi_{i}(\varepsilon)\varphi_{k-1-i}(\varepsilon)d\varepsilon,
\end{eqnarray}
in which  $Q_{k-1}$ is determined by the restriction condition (\ref{restrict}).   Consequently, $\varphi_{k}(y)$ and $S_{k}(y)$ of Eqs.~(\ref{high:1}) and (\ref{high:2}) are obtained.  Then, we have the $n$th-order approximation:
\begin{eqnarray}
  \varphi(y) &=& \sum_{k=0}^{n}\varphi_{k}(y),\label{w:sum:phi:n}\\
   S(y) &=& \sum_{k=0}^{n}S_{k}(y),\label{w:sum:s:n}\\
   Q &=& \sum_{k=0}^{n}Q_{k}.\label{w:sum:Q:n}
\end{eqnarray}

According to \cite{Liaobook}, the initial guesses $\varphi_{0}(y)$ and $S_{0}(y)$ should obey the solution expression (\ref{homotopy:series}) and satisfy the boundary conditions. Thus, we choose
\begin{eqnarray}
\varphi_{0}(y)&=&\frac{-2a}{2\lambda+1}[(\lambda+1)y-y^{2}],\\
S_{0}(y)&=&0,  \label{w:initial:s}
\end{eqnarray}
as the initial guesses of $\varphi(y)$ and $S(y)$, respectively.

\subsection{Results given by the HAM-approach without iteration}

\begin{figure}
    \begin{center}
        \begin{tabular}{cc}
            \includegraphics[width=3in]{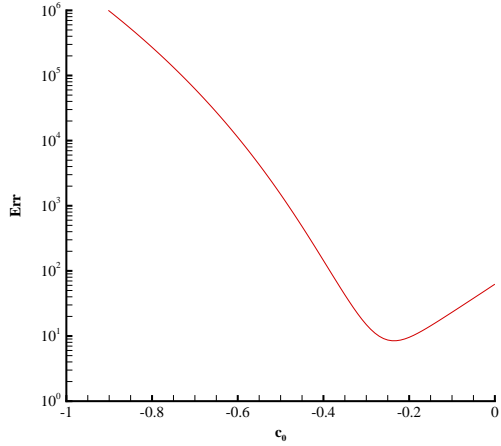} \\
        \end{tabular}
    \caption{The sum of the squared residual errors $Err$ versus $c_{0}$ in the case of $a = 5$ for a circular plate with the clamped boundary, gained by the HAM-based approach (without iteration).} \label{HAM:a5:opt:c0:5th}
    \end{center}
\end{figure}

First, the HAM-based approach (without iteration) for given central deflection is used to solve the Von K{\'a}rm{\'a}n's plate equations in the integral form with the clamped boundary.  Without loss of generality, let us consider the same case of $a=5$, equivalent to the central deflection $w(0)/h=3$. As shown in Fig.~\ref{HAM:a5:opt:c0:5th}, the sum of the squared residual errors $Err$ arrives its minimum at $c_{0}\approx-0.25$.  As shown in Table~\ref{Table:HAM:noniter:a5}, the sum of the squared residual errors quickly decreases to $3.6\times10^{-7}$  by means of $c_{0}\approx-0.25$ in the case of $a=5$.  Note that, the Chien's perturbation method \cite{Qian} (using $a$ as the perturbation quantity) is only valid for $w(0)/h<2.44$, equivalent to $a<4$, for a plate with the clamped boundary. This  again illustrates the superiority of the HAM-based approach over the perturbation method.

\begin{table}
\tabcolsep 0pt
\caption{The sum of the squared residual errors $Err$ and the corresponding load $Q$ versus the  order of approximation  in the case of $a= 5$ for a circular plate with the clamped boundary, given by the HAM without iteration using $c_{0} = -0.25$.}
\vspace*{-12pt}\label{Table:HAM:noniter:a5}
\begin{center}
\def\temptablewidth{0.8\textwidth}
{\rule{\temptablewidth}{1pt}}
\begin{tabular*}{\temptablewidth}{@{\extracolsep{\fill}}cccccc}
Order of approx. &  $Err$    &  $Q$                    \\
\hline
20 &  $3.4\times10^{-2}$  &  132.3                          \\
40 &  $2.2\times10^{-3}$  &  132.5                  \\
60 &  $8.8\times10^{-5}$ &  132.3                  \\
80 &  $9.7\times10^{-7}$ &  132.2                         \\
100 & $3.6\times10^{-7}$   &  132.2                      \\
\end{tabular*}
{\rule{\temptablewidth}{1pt}}
 \end{center}
 \end{table}

Similarly, for given value of $a$, we can always find its optimal value of $c_{0}$, which can be expressed by the empirical formula:
\begin{equation}
c_{0}=-\frac{11}{11+a^2}~~~~~~~(a\leq5).\label{c0:iter:W}
\end{equation}
In addition, the convergent results of the external load $Q$ in case of different values of $a$ for a circular plate with the clamped boundary are given in Table~{\ref{HAM:noniter:diff:a}}.

\begin{table}
\tabcolsep 0pt
\caption{The convergent results of the external load $Q$ in case of different values of $a$ for a circular plate with the clamped boundary, given by the HAM-based approach without iteration using the optimal convergence-control parameter $c_{0}$ given by (\ref{c0:iter:W}).}
\vspace*{-12pt}\label{HAM:noniter:diff:a}
\begin{center}
\def\temptablewidth{0.8\textwidth}
{\rule{\temptablewidth}{1pt}}
\begin{tabular*}{\temptablewidth}{@{\extracolsep{\fill}}cccccc}
$a$  & $c_{0}$  &  $Q$       \\
\hline
1    & -0.92   &  4.8       \\
2   & -0.73    & 14.6      \\
3   & -0.55    & 35.2    \\
4   & -0.41   &  72.4     \\
5   & -0.31   &   132.2   \\
\end{tabular*}
{\rule{\temptablewidth}{1pt}}
 \end{center}
 \end{table}

\subsection{Convergence acceleration by means of iteration}

As shown in \S~2.3,   iteration can greatly accelerate the convergence of the homotopy-series solutions.  Therefore, we used the HAM-based iteration approach for given central deflection to solve the Von K{\'a}rm{\'a}n's plate equations with the clamped boundary in this subsection.   The iteration order $M$ and the truncation order $N$ are defined in a similar way to Eqs.~(\ref{Q:M:order:phi}) and (\ref{Q:nterms2}).

\begin{figure}[t]
    \begin{center}
        \begin{tabular}{cc}
            \includegraphics[width=2.5in]{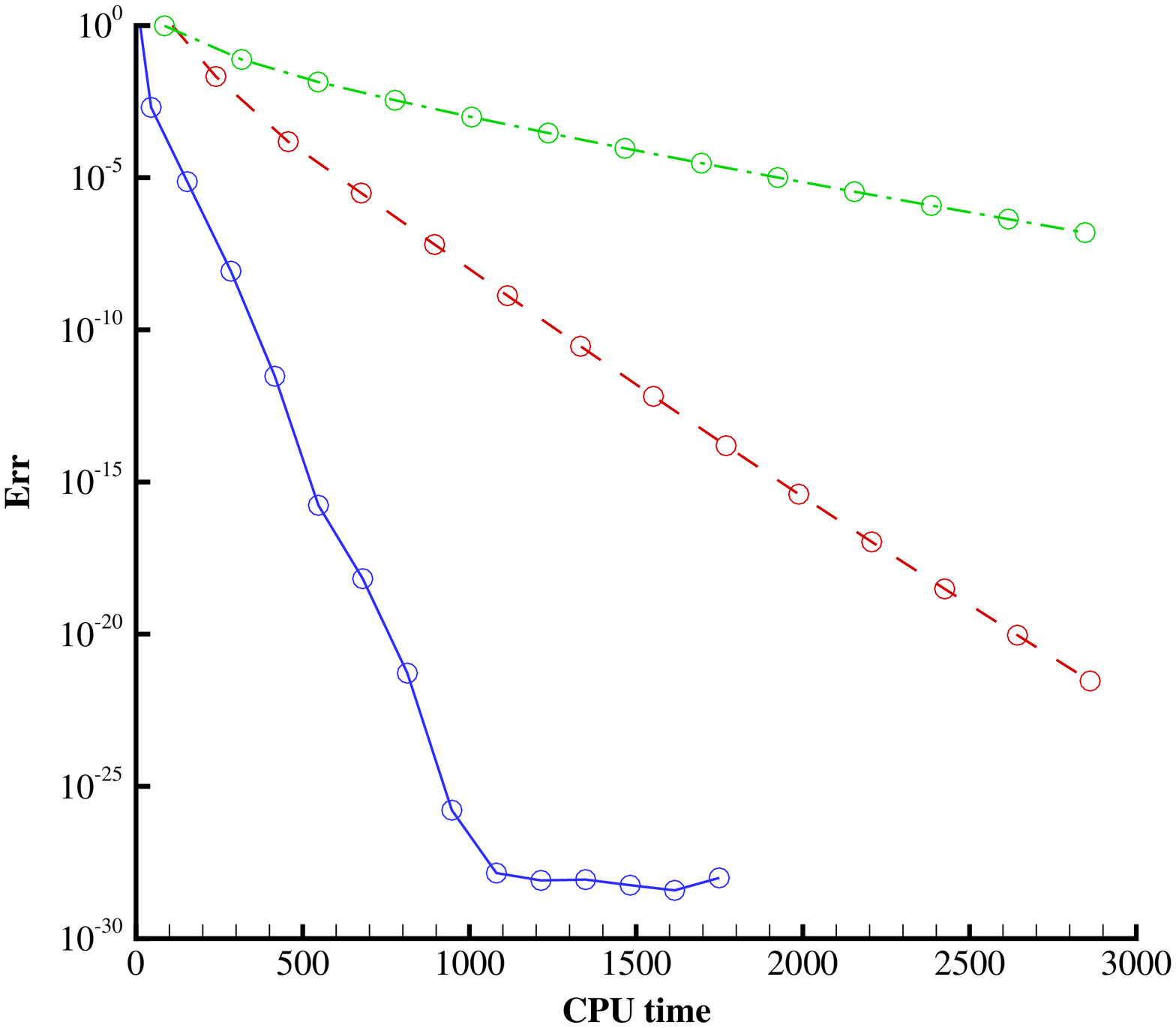} \\
        \end{tabular}
    \caption{The sum of the squared residual errors $Err$ versus the CPU time in the case of $a = 5$ for a circular plate with the clamped boundary, gained by the interpolation iterative method \cite{Keller, Zheng}, and the HAM-based approach for given external load $Q$ and central deflection $a$, respectively.  Dash-dotted line: results given by the interpolation iterative method \cite{Keller, Zheng} using the interpolation parameter $\theta=0.1$; Dashed line: results given by the HAM-based 5th-order iteration approach for given external load $Q$ using  the optimal convergence-control parameter $c_{0}=-0.15$ and the truncation order $N=100$; Solid line: results given by the HAM-based 5th-order iteration approach for given central deflection $a$ using the optimal convergence-control parameter $c_{0}=-0.5$ and the truncation order $N=100$.} \label{W:versus:Q:CPU}
    \end{center}
\end{figure}

\begin{figure}
    \begin{center}
        \begin{tabular}{cc}
            \includegraphics[width=2.5in]{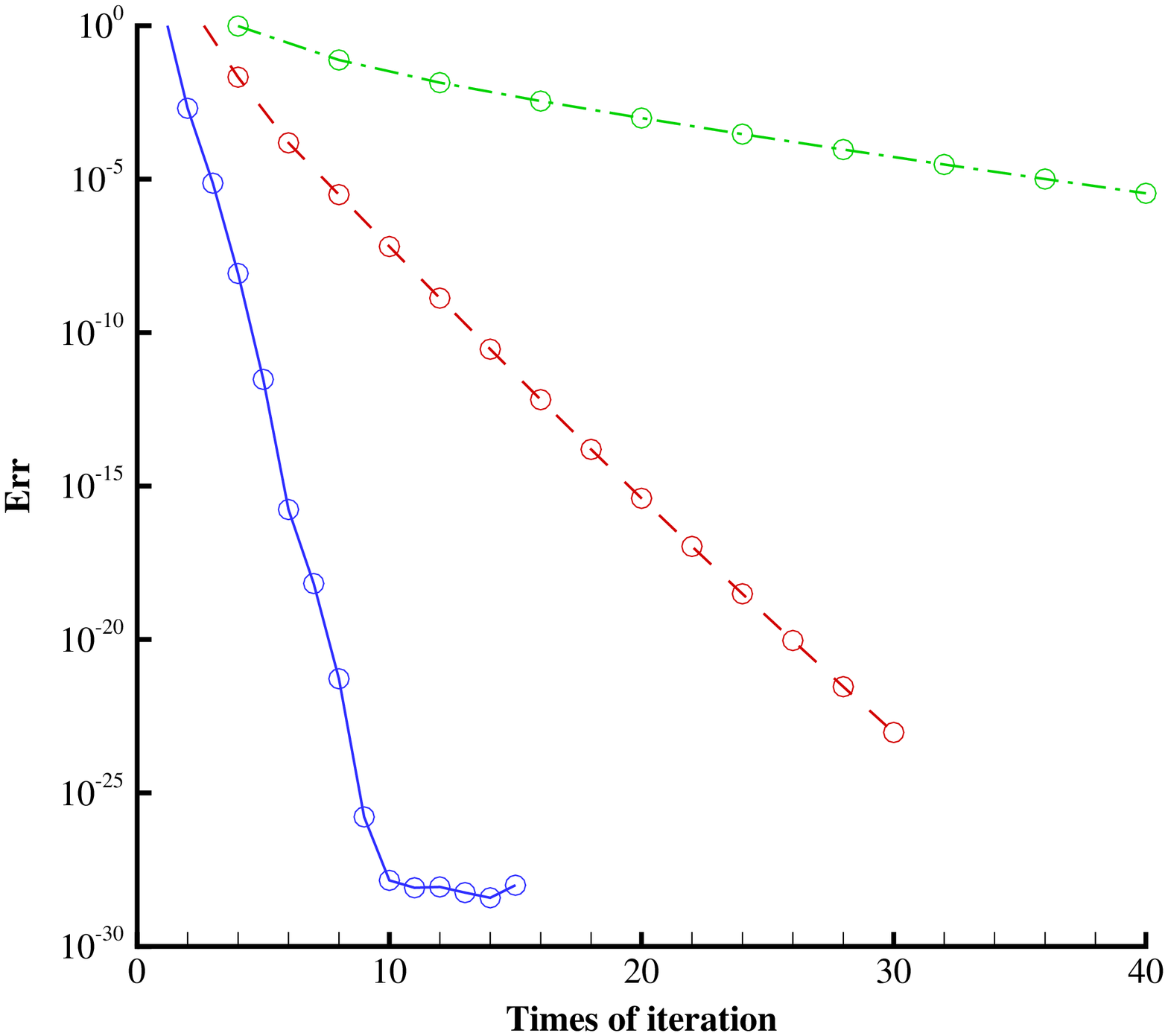} \\
        \end{tabular}
    \caption{The sum of the squared residual errors $Err$ versus iteration times in the case of $a = 5$ for a circular plate with the clamped boundary, gained by the interpolation iterative method \cite{Keller, Zheng},  and the HAM-based approach for given external load $Q$ and central deflection $a$, respectively. Dash-dotted line: results given by the interpolation iterative method \cite{Keller, Zheng} using the interpolation parameter $\theta=0.1$; Dashed line: results given by the HAM-based 5th-order iteration approach for given external load $Q$ using the optimal convergence-control parameter $c_{0}=-0.15$ and the truncation order $N=100$; Solid line: results given by the HAM-based 5th-order iteration approach for given central deflection $a$ using the optimal convergence-control parameter $c_{0}=-0.5$ and the truncation order $N=100$.} \label{W:versus:Q:order}
    \end{center}
\end{figure}

The iteration order $M=5$ and the truncation order $N=100$ are used in all cases described below.  Without loss of generality, let us first consider the same case of $a=5$. Note that the sum of the squared residual errors quickly decreases to $1.4\times10^{-28}$ in only $10$ iteration times,  as shown in Table~\ref{Table:HAM:iter:a5}.  In addition, as shown in Figs.~\ref{W:versus:Q:CPU} and \ref{W:versus:Q:order},  the convergence is much faster by introducing $W(0)$ into the Von K{\'a}rm{\'a}n's plate equations.   Therefore, using the central deflection  indeed makes sense.  It should be emphasized that the HAM-based iteration approaches converges much faster than the interpolation iterative method \cite{Keller, Zheng}, as shown in Figs.~\ref{W:versus:Q:CPU} and \ref{W:versus:Q:order}.  

\begin{table}
\tabcolsep 0pt
\caption{The sum of the squared residual $Err$ and the load $Q$ versus the iteration times in the case of $a= 5$ for a circular plate with the clamped boundary, given by the HAM-based 5th-order iteration approach using $c_{0} = -0.5$ and the truncation order $N=100$.}
\vspace*{-12pt}\label{Table:HAM:iter:a5}
\begin{center}
\def\temptablewidth{0.8\textwidth}
{\rule{\temptablewidth}{1pt}}
\begin{tabular*}{\temptablewidth}{@{\extracolsep{\fill}}cccccc}
$m$, times of iteration. &  $Err$    &  $Q$                    \\
\hline
2 &  $2.0\times10^{-3}$  &  132.0                          \\
4 &  $8.4\times10^{-9}$  &  132.2                  \\
6 &  $1.7\times10^{-16}$ &  132.2                  \\
8 &  $5.3\times10^{-22}$ &  132.2                         \\
10 & $1.4\times10^{-28}$   &  132.2                      \\
\end{tabular*}
{\rule{\temptablewidth}{1pt}}
 \end{center}
 \end{table}

\begin{figure}[t]
    \begin{center}
        \begin{tabular}{cc}
            \includegraphics[width=3in]{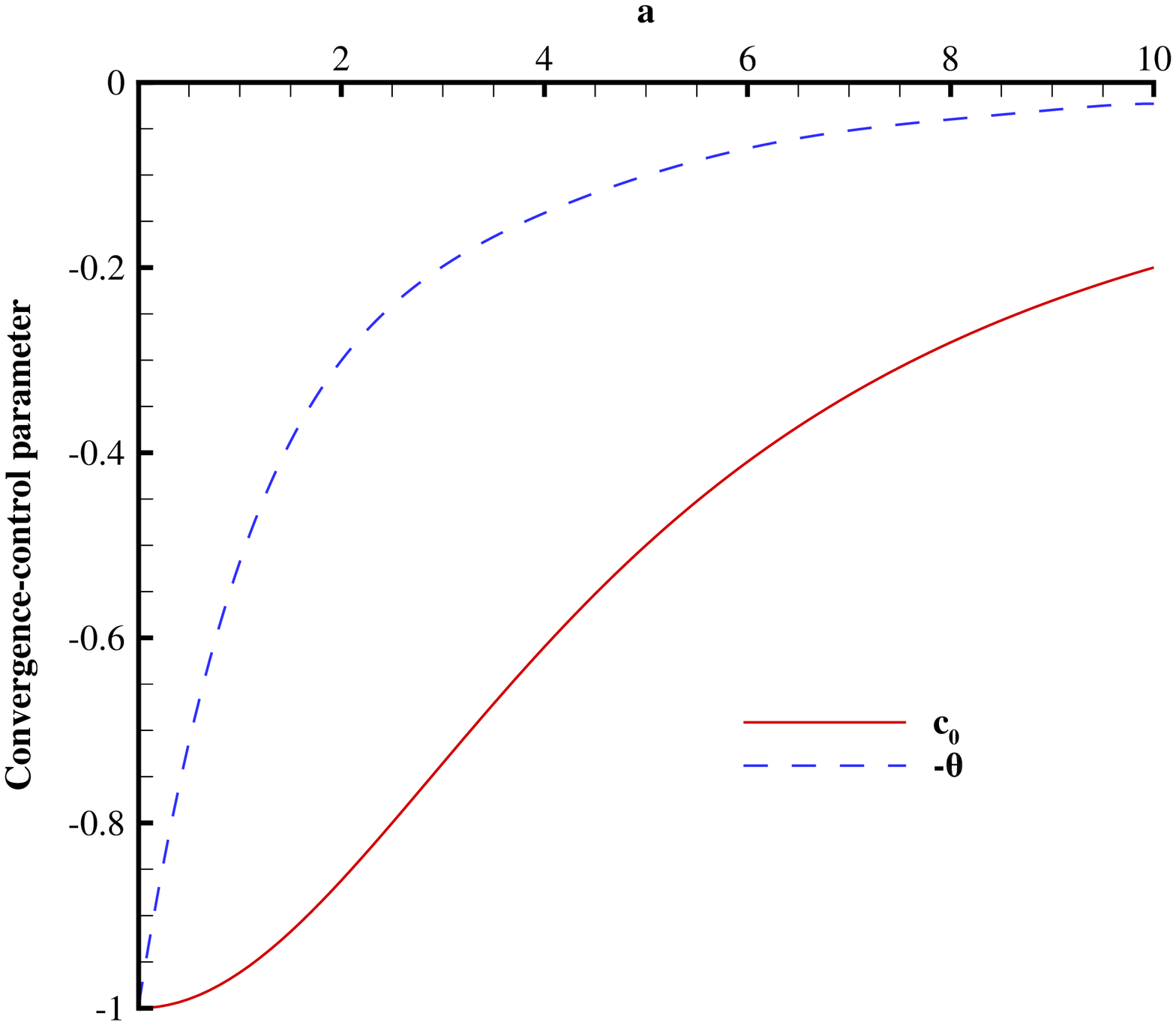} \\
        \end{tabular}
    \caption{The convergence-control parameter $c_{0}$ and the interpolation parameter $\theta$ versus the central deflection $a$ for a circular plate with the clamped boundary, gained by the interpolation iterative method \cite{Zheng} for given external load $Q$ and the HAM-based 5th-order iteration approach for given central deflection $a$ using the truncation order $N=100$. Dashed line: $-\theta$; Solid line: $c_{0}$.} \label{curve:c}
    \end{center}
\end{figure}

Similarly, for any a given value of $a$, we can always find the optimal value of the convergence-control parameter $c_{0}$,  which can be expressed by the empirical formula:
\begin{equation}
c_{0}=-\frac{25}{25+a^2}. \label{c0:iter:w}
\end{equation}
 As shown in Table~\ref{HAM:iter:diff:a}, convergent solutions can be obtained even in the case of $a=30$, equivalent to $w(0)/h=18.2$.   As shown in Fig.~\ref{curve:c},  as $a$ increases,   the interpolation parameter $\theta$ of the interpolation iterative method \cite{Keller, Zheng} approaches to $0$ much faster than the optimal convergence-control parameter $c_{0}$ mentioned-above.  This explains why the HAM-based iteration approaches converge much faster than the interpolation iterative method \cite{Keller, Zheng}.  

\begin{table}
\tabcolsep 0pt
\caption{The convergent homotopy-approximation of $Q$ in case of different values of $a$ for a circular plate with the clamped boundary, given by the HAM-based 5th-order iteration approach with the optimal convergence-control parameter $c_{0}$ given by (\ref{c0:iter:w}) and the truncation order $N=100$.}
\vspace*{-12pt}\label{HAM:iter:diff:a}
\begin{center}
\def\temptablewidth{0.8\textwidth}
{\rule{\temptablewidth}{1pt}}
\begin{tabular*}{\temptablewidth}{@{\extracolsep{\fill}}cccccc}
$a$  & $c_{0}$  &  $Q$       \\
\hline
5    & -0.50   &  132.2        \\
10   & -0.20    &  957.7      \\
15   & -0.10    &  3152.1     \\
20   & -0.06   &   7386.9     \\
25   & -0.04   &   14334.1    \\
30   & -0.03   &   24665.7     \\
\end{tabular*}
{\rule{\temptablewidth}{1pt}}
 \end{center}
 \end{table}

\section{Concluding remarks}  

In this paper, the homotopy analysis method (HAM) is successfully applied to solve the Von K{\'a}rm{\'a}n's plate equations in the integral form for a circular plate  with the clamped boundary under an arbitrary uniform external pressure.   Two HAM-based approached are proposed.  One is for given load $Q$, the other for given central deflection.  Both of them are valid for arbitrary uniform external pressure, by mens of choosing a proper value of the so-called convergence-control parameters $c_1$ and $c_2$ in the frame of the HAM.  Besides,  it is found that iteration can greatly accelerate the convergence of solution series.   In addition,  it is proved that the interpolation iterative method \cite{Keller, Zheng} is a special case of the HAM-based $1$st-order iteration approach for given external load $Q$ when $c_{1}=-\theta$ and $c_{2}=-1$, where $\theta$ denotes the interpolation parameter of the  interpolation iterative method\footnote{Note that,  it  was proved in Part (I) that the perturbation methods for {\em any} a perturbation quantity (including Vincent's \cite{Vincent} and Chien's \cite{Qian} perturbation methods) and the modified iteration method \cite{YehLiu} are only the special cases of the HAM-based approaches  when $c_{0}=-1$,  for the Von K{\'a}rm{\'a}n's plate equations in the differential form for a circular plate under a uniform external pressure.}.   Therefore,  like Zheng and Zhou \cite{Zheng2},  one can similarly prove  that  the HAM-based approaches for the Von K{\'a}rm{\'a}n's plate equations in the integral form are valid for an arbitrary uniform external pressure, at least in some special cases such as $c_{1}=-\theta$ and $c_{2}=-1$.      Finally, it should be emphasized that  the HAM-based iteration approaches converge much faster than the interpolation iterative method \cite{Keller, Zheng}, as shown in Figs.~\ref{W:versus:Q:CPU} and \ref{W:versus:Q:order}.  All of these illustrate the validity and potential of the HAM for the famous Von K{\'a}rm{\'a}n's plate equations, and show the superiority of the HAM over perturbation methods.    Without doubt, the HAM can be  applied to solve some challenging  nonlinear problems in solid mechanics.

\section*{Acknowledgment}
This work is partly supported by National Natural Science Foundation of China (Approval No. 11272209 and 11432009) and State Key Laboratory of Ocean Engineering (Approval No. GKZD010063).

\section*{References}
\bibliographystyle{elsarticle-num}
\bibliography{plate}
\end{document}